\documentclass[11pt]{amsart}
\usepackage{mathrsfs}
\usepackage{xcolor}

\usepackage{amsmath,amssymb,amsthm,amsfonts,graphicx,color}
\usepackage[font=small,labelfont={bf,sf}, textfont={sf},margin=0em]{caption}
\usepackage{url,hyperref}
\definecolor{limegreen}{rgb}{0.196,0.804,0.196}
\definecolor{darkgreen}{rgb}{0.0,0.5,0.0}
\definecolor{darkbluegreen}{rgb}{0,0.3,0.6}
\definecolor{badgerred}{rgb}{0.715,0.004,0.004}
\hypersetup{colorlinks, 
citecolor=badgerred,
filecolor=black,
linkcolor=blue,
urlcolor=darkgreen,
bookmarksopen=false,
pdftex}

\newcommand{\cM}{\mathcal{M}}

\newcommand{\R}{{\mathbb R}}

\newcommand{\calM}{\mathcal{M}}
\newcommand{\calC}{\mathcal{C}}
\newcommand{\bx}{\mathbf{x}}
\newcommand{\by}{\mathbf{y}}
\newcommand{\bO}{\mathbf{0}}
\newcommand{\barr}{\bar}

\newtheorem{theorem}{Theorem}[section]
\newtheorem{lemma}[theorem]{Lemma}
\newtheorem{proposition}[theorem]{Proposition}
\newtheorem{definition}[theorem]{Definition}

\newtheorem{remark}[theorem]{Remark}

\newtheorem{claim}[theorem]{Claim}

\numberwithin{equation}{section}
\numberwithin{theorem}{section}

\begin{document}
\title[Stability of neckpinch singularities]{Stability of neckpinch singularities}
\author{Felix Schulze and Natasa Sesum}
\address{\noindent FS: Department of Mathematics, Zeeman Building, University of Warwick, Gibbet Hill Road, Coventry CV4 7AL, UK}
\address{\noindent NS: Department of Mathematics, Rutgers University, Piscataway NJ 08854, USA}

\begin{abstract}
In this paper, we study the stability of neckpinch singularities. We show that if a mean curvature flow $\{M_t\}$ develops only finitely many neckpinch singularities at the first singular time, then the mean curvature flow starting at any sufficiently small perturbation of $M_0$ can also develop only neckpinch type singularities at the first singular time. We also show stability of nondegenerate neckpinch singularities in the above sense, which speaks in favor of stability of Type I singularities. 
\end{abstract}
\thanks{The first author was supported by a Leverhulme Trust Research Project Grant RPG-2016-174. The second author was supported by the National Science Foundation under grants DMS-1056387 and DMS-1811833.} 
\maketitle 

\section{Introduction}

Let $F : M \times [0,T) \to \mathbb{R}^{n+1}$ denote a family of embeddings of closed hypersurfaces evolving by mean curvature flow
\[\left(\frac{\partial}{\partial t} F\right)^\perp = \vec{H} =  - H\, \nu\]
and let  $M_t = F(M^n,t)$. 

It is well known that the mean curvature flow starting at a closed hypersurface becomes singular in finite time. One of the most important problems in the mean curvature flow theory is understanding singularities as they arise. Huisken conjectured that a generic mean curvature flow has only spherical and cylindrical singularities. The breakthrough work towards proving Huisken's conjecture was made by Colding and Minicozzi in \cite{CM12} who proved that spheres and cylinders are the only
linearly stable singularity models for mean curvature flow. In particular, all other
singularity models are linearly unstable and are expected to occur only non-generically. In recent work of the second author together with Chodosh, Choi and Mantoulidis \cite{CCMS}, they make a significant contribution towards verifying Huisken's conjecture in the case of surfaces in $\mathbb{R}^3$. More precisely they show   that a large class of unstable singularity models are, in fact,
avoidable by a slight perturbation of the initial data. Roughly stated, they show that
the mean curvature flow of a generic closed embedded surface in $\mathbb{R}^
3$
encounters only spherical and cylindrical singularities until the first time
it encounters a singularity modelled on a self-shrinker with multiplicity $\ge 2$, or a self-shrinker which has a
cylindrical end but which is not globally a cylinder. It is conjectured that none of those last two scenarios actually happen.

Assume the flow $\{M_t\}$ has a singularity at $(x_0,T)$, where $T$ is the first singular time. Denote its spacetime track by $\mathcal{M} = \cup_{t \in [0,T)} M_t\times \{t\} \subset \R^{n+1}\times \R$.

Let us explain what we mean by having a limit flow at $(x_0,T)$. For any $X = (x,t) \in \mathcal{M}$ and $\lambda > 0$ we denote by $\calM_{X,\lambda} = \mathcal{D}_\lambda(\calM - X)$ the flow which is obtained from $\mathcal{M}$ by translating $X$ to the origin and parabolically dilating by $\lambda$. Given $X\in \mathcal{M}$, and arbitrary sequences $X_i = (x_i,t_i)\to X$, and $\lambda_i\to\infty$, one can always pass to a (weak) subsequential limit of $\mathcal{M}_{X_i,\lambda_i}$. We will call any such limit $\mathcal{M}_{\infty}$ a \emph{limit flow} at $X$. In case that for all $i$ it holds $t_i \leq T$ with $t_i<T$ if $x_i \neq x$ we call $\mathcal{M}_{\infty}$ a \emph{special limit flow} at $X$. In the  special case when $X_i = X$ for all $i$, the blow up limit $\mathcal{M}_{\infty}$ is called a \emph{tangent flow}. This is nearly identical with White's definition, see \cite{White03}. 

We will in the following denote the space time track of the self-similar evolution of a unit multiplicity cylinder $S^{n-1}(\sqrt{2(n-1)})\times \mathbb{R}$ with axis $a \in \mathbb{RP}^{n+1}$ as 
$$ \calC = \cup_{t<0} (S^{n-1}(\sqrt{2(n-1)(-t)})\times \mathbb{R})\times \{t\} \subset \R^{n+1}\times \R\, .$$
\begin{definition}
\label{def-nondegenerate}
Let $\mathcal{M}$ be a unit regular, intergal Brakke flow. We say $\mathcal{M}$ has a neckpinch singularity at $X$ if a tangent flow at $X$  is equal to $\calC$. We say it is a {\bf nondegenerate} neckpinch if every nontrivial special limit flow at $X$ is up to rotation and translation in space-time equal to $\calC$ with unit multiplicity. In the case at least one of the nontrivial special limit flows at $X$ is not a round cylinder, we say it develops a {\bf degenerate} neckpinch at $X$.
\end{definition}

\begin{remark} We point that our notion of a nondegenerate neckpinch singularity does not fully agree with the more analytical notion previously employed for example in \cite{AV,GS}. Nevertheless, we feel that our notion reflects what is geometrically understood as a generate neckpinch and thus correspondingly as a nondegenerate neckpinch.
 
\end{remark}
Our goal in this work is to show that "nondegenerate neckpinch" type singularities are stable. We will first show that local neckpinch singularities are stable in the sense of the theorem below.

\begin{theorem}
\label{thm-main1}
Let $\{M_t\}_{t\in [0,T)}$ be a smooth mean curvature flow that develops at most finitely many isolated neckpinch singularities at time $T$.  Assume each neck locally disconnects the manifold into two pieces none of which disappears at time $T$. There exists a $\delta > 0$, so that for every embedding $\barr{F}_0: M\to \mathbb{R}^{n+1}$, for which $\|\barr{F}_0(\cdot) - F(\cdot,0)\|_{C^2} < \delta$, the tangent flows at the first singular time of the mean curvature flow $\barr{M}_t$, starting at $\barr{M}_0 = \barr{F}(M)$ are only  multiplicity one round cylinders $\calC$.
\end{theorem}

\smallskip 

The next result states that local nondegenerate neckpinch singularities are stable, which is not expected to be the case for degenerate neckpinches.

\begin{theorem}
\label{thm-main2}
Assume $\{M_t\}_{t\in [0,T)}$ is a smooth mean curvature flow that develops at most finitely many isolated  nondegenerate neckpinch singularities at time $T$. Then, there exists a $\delta > 0$, so that for every embedding $\barr{F}_0: M\to \mathbb{R}^{n+1}$, for which $\|\barr{F}_0(\cdot) - F(\cdot,0)\|_{C^2} < \delta$, the mean curvature flow $\barr{M}_t$, starting at $\barr{M}_0 = \barr{F}(M)$ can develop only nondegenerate neckpinch singularities at the first singular time. 
\end{theorem}

To prove Theorem \ref{thm-main2} we use Theorem \ref{thm-main1}, once we show the following result to hold for isolated nondegenerate neckpinches.

\begin{theorem}
\label{thm-no-vanishing}
Assume $\{M_t\}_{t\in [0,T)}$ is a smooth mean curvature flow that develops at most finitely many isolated  nondegenerate neckpinch singularities at time $T$. Then each of the necks locally disconnects the manifold into two pieces none of which disappears at time $T$. 
\end{theorem}

The proofs of Theorems \ref{thm-main1}, \ref{thm-main2} and \ref{thm-no-vanishing} heavily rely on the proof of the mean convex neighborhood theorem by Choi, Haslhofer, Hershkovits and White \cite{CHHW}, see also \cite{CHH18}. 

\begin{remark}
 Both Theorem \ref{thm-main1} and Theorem \ref{thm-main2} are direct consequences of corresponding local statements which do not depend on perturbations of the initial condition, see Theorem \ref{thm-main1-local}  and Theorem \ref{thm-nondegenerate-stability-local}.
\end{remark}

The stability of neckpinch singularities for mean convex mean curvature flow has been investigated by White \cite{White13}, where he shows (among other results) that there is an open set of mean convex initial conditions developing cylindrical neckpinches of any generalised cylindrical type. 

It is helpful to clarify the relation between Type I, resp.~Type II, singularities and nondegenerate and resp. degenerate neckpinches. We recall the definition of a Type I, resp.~Type II, singularity.

\begin{definition}
\label{def-typeI}
A mean curvature flow $\mathcal{M}$ has a (local) type I singularity at point $X_0=(\mathbf{x}_0,t_0)$ if $X_0$ is a singular point of the flow and
there exists $\delta>0$ and $C>0$ such that
$$|A|^2(\mathbf{x},t) \leq \frac{C}{t_0-t}$$
for all $(x,t) \in B_\delta(\mathbf{x}_0)\times (t_0-\delta^2,t_0)$.
Otherwise we call $X_0$ a (local) type II singularity.
\end{definition}

We have the following theorem.

\begin{theorem}\label{prop-typeI}
Assume that a smooth mean curvature flow $(M_t)_{t\in [0,T)}$ has a neckpinch singularity at $X=(x_0,T)$. Then $X$ is a (local) type I singularity if $X$ is a nondegenerate neckpinch. Similarly $X$ is a (local) type II singularity if $X$ is a degenerate neckpinch.
\end{theorem}

The theorem follows directly from Hamilton's rescaling procedure and \cite{CHHW}, we give the proof of it in Section \ref{sec-typeI}. 
From that point of view, we can see Theorem \ref{thm-main2} speaking in favor of  stability of  Type I singularities. More precisely, it is expected that if the flow $\calM$ develops only Type I singularities at the first singular time $T$, then the flow starting at any sufficiently small perturbation of $M_0$ can also develop only Type I singularities.

In \cite{GS} the authors constructed an open set of closed rotationally symmetric solutions to the mean curvature flow that develop a nondegenerate neckpinch singularity in finite time $T < \infty$. More precisely, they constructed rotationally symmetric solutions defined by a map $r = U(x,t)$, where
\[\frac{\partial}{\partial t} U = \frac{U_{xx}}{1+u_x^2} - \frac{n-1}{U}\]
so that, the rescaling 
$$u(y,\tau) := \frac{U(x,t)}{\sqrt{T-t}}, \qquad y= \frac{x}{\sqrt{T-t}}, \,\, \tau = - \log (T-t)$$ 
satisfies,  as $\tau\to \infty$,  the asymptotics 
\begin{equation}
\label{eq-profile}
u(y,\tau) = \sqrt{2(n-1)} + \frac{\sqrt{2(n-1)}}{4\tau} \, (y^2 - 2) + o(|\tau|^{-1}) \quad \mbox{on}\,\, |y| \le L
\end{equation}
for every number  $0 < L < \infty$. These are some of many examples of mean curvature flows to which the theorems above could be applied. 

\begin{remark}
With no loss of generality we will assume in the proofs below that our flow has only one neckpinch singularity at $(x_0,t_0)$. The analysis is the same in the case we have finitely many of them.
\end{remark}

\medskip 

The organization of the paper is as follows. In section \ref{sec-rot-symm}, using the fact that our flow has isolated neckpinches at the first singular time, we show several results to hold on a macroscopic scale for a perturbed flow. We end the section with the proof of Theorem \ref{thm-no-vanishing}. In section \ref{sec-neck} we prove Theorem \ref{thm-main1}, while Theorem \ref{thm-main2} we show in section \ref{sec-nondeg-neck}.  In section \ref{sec-typeI} we give a proof of Theorem \ref{prop-typeI}. 

\medskip

{\bf Acknowledgements:} The authors would like to thank Panagiota Das\-ka\-lopou\-los for many insightful and helpful discussions.

\section{Almost rotationally symmetric neighborhoods of cylindrical singularities} 
\label{sec-rot-symm}

Since we will give local statements for our perturbative results, we will in this section assume that the flow $\calM$ is a unit-regular, integral Brakke flow which has at the point $({\bf 0},0)$ a multiplicity one cylinder $\calC$ 
as its unique tangent flow. For background on weak mean curvature flows we refer the reader to \cite[Section 2]{CCMS} and \cite[Section 2]{CHH18}.

Recall the definition of a mean convex mean curvature flow being $\alpha$-noncollapsed.
\begin{definition}\label{def-alpha-part}
A smooth mean convex hypersurface $M$ bounding an open region $\Omega$ in $\mathbb{R}^{n+1}$ is said to be $\alpha$-noncollapsed (on the scale of its mean curvature) if for every $x\in M$ there are balls $B_{int}\subset \Omega$ and $B_{ext} \subset \mathbb{R}^{n+1}\backslash Int(K)$ of radius at least $\alpha/H(x)$ that are tangent to $M$ at $x$ from the interior and exterior of $\Omega$, with $x\in \partial\Omega = M$. We say that a smooth mean curvature flow $M_t$ is $\alpha$-noncollapsed  if every one of its time slices is $\alpha$-noncollapsed. We say that a weak mean curvature flow is $\alpha$-noncollapsed if all of its regular points have positive mean curvature and the support of its time-slices satisfy the above condition at every regular point.
\end{definition}

\smallskip 

The proofs of Theorems \ref{thm-main1} and \ref{thm-main2} heavily use  the following Proposition whose proof relies on the proof of the mean convexity neighborhood theorem, see \cite{CHHW}. We use the same terminology as in the statement there. We recall that by \cite{CM} multiplicity one cylindrical tangent flows are unique. 

We also denote with $R(X)$ the regularity scale of a flow $\calM$ at $X$, defined as in \cite{CHH18, CHHW}:

\begin{definition}
\label{def-reg-scale}
Let $S\subset \mathbb{R}^{n+1}$ and $p\in S$. Define the {\it regularity scale} $R(S,p)$ to be the supremum of $r > 0$ so that $S\cap B(p,r)$ is a smooth $n$-dimensional manifold {\rm(}with no boundary in $B(p,r)${\rm)} properly embedded in $B(p,r)$ and such that the norm of the second fundamental form at each point of $S\cap B(p,r)$ is $\le 1/r$. If there is no such $r$, we let $R(S,p) = 0$. If $X=(\bx,t)$ is a regular point of a mean curvature flow $\calM$ we let $R(X)$ be the supremum of $r>0$ such that $R(M_t, \bx) \leq r$ for all $t \in (t-r^2,t)$.
\end{definition}

We have the following structural statement about flows close to a flow with a neckpinch singularity.

\begin{proposition}
\label{prop-help}
Let $\calM$ be a unit-regular, $n$-integral Brakke flow which has at $O=(\bO,0)$ a multiplicity one cylinder $\calC$ as its unique tangent flow. If $n=2$ we also assume that $\calM$ is cyclic mod 2. Let $\calM^i$ be a sequence of unit-regular, integral Brakke flows converging locally around the point $O$
to $\calM$. Then there exist $r, \eta, \alpha>0$, $i_0 \in \mathbb{N}$, such that for all $i \geq i_0$
\begin{enumerate}
\item[$ (i)$] $H \ge \eta$ for every regular point of $\calM^i$ in $\bar{B}(\bO, r) \times [-r^2, r^2]$.
\item[$(ii)$]The flows $\calM^i$ have only multiplicity one neck and spherical singularities in $\bar{B}(\bO, r) \times [-r^2, r^2]$.
\item[$(iii)$] The flows $\calM^i$ are smooth in $\bar{B}(\bO, r)$ for a.e.~$t\in [-r^2, r^2]$ and are smooth outside a set of Hausdorff dimension 1 for every $t\in [-r^2, r^2]$.
\item[$(iv)$] If $X_i \in \calM^i$ are regular points and $X_i \rightarrow O$, then any subsequential limit of $\calM^i_{X_i, R(X_i)^{-1}}$ is either a round shrinking sphere, a round shrinking cylinder, a translating bowl, or an ancient oval, with the axis aligned with the axis of the neck singularity  of $\cM$ at $O$.
\item[$(v)$]  $\lambda_1 + \lambda_2 \ge  \eta H$ for every regular point of $\calM^i$ in $\bar{B}(\bO, r) \times [-r^2, r^2]$.
\item[$(vi)$] The flows $\calM^i$ are $\alpha$-noncollapsed in $\bar{B}(\bO, r) \times [-r^2, r^2]$.
\end{enumerate}
\end{proposition}

\smallskip

\begin{remark} Note that the weak convergence of Brakke flows is metrizable, see for example \cite{SW}. Thus it is possible to quantify this in the above statement as well.
\end{remark}

\begin{proof} [Proof of Proposition \ref{prop-help}]
The proof is a direct extension of the proof of  \cite[Proposition 9.1]{CHHW}. We explain how to slightly change the set-up to allow for the above conclusions. 

Recall the definition of the cylindrical scale $Z(X)$ at a point $X=(\bx,t) \in \mathcal{M}$ of an asymptotically cylindrical ancient flow,
 \cite[Definition 4.11]{CHHW}. In our case we fix in the definition the axis $a$ of the asymptotic cylinder $\mathcal{C}$, which is the axis of the neck singularity of $\calM$ at $O$. Slightly informally speaking, compare \cite[Proof of Proposition 2.1]{CHH19}, one defines $Z(X)$ as follows:  Fix $\varepsilon>0$ sufficiently small, and consider the rescaled flow $\mathcal{M}_{X,1/r_\alpha}$ on dyadic annuli of radius $r_\alpha = 2^\alpha$ for $\alpha \in \mathbb{Z}$. Then $Z(X)$ is the first cylindrical scale, i.e.~the infimum of $r_\alpha$ such that $\mathcal{M}_{X,1/r_\alpha}$ is $\varepsilon$-close to $\mathcal{C}$.
 
Recall that the tangent flow of $\calM$ at $O$ is the unique multiplicity one cylinder $\mathcal{C}$. Thus for every $r>1$ there exists $\Lambda(r) >0$ such that for all $\lambda \geq \Lambda(r)$ we have that
$$ \calM_{O,\lambda} \text{ is } (2r)^{-1} \text{ close in } C^{2,\alpha} \text{ to } \mathcal{C} \text{ on } B_{2r}(\bO)\times \left(-(2r)^{2}, - (2r)^{-2}\right).$$

Since the $\calM^i$ are converging weakly to $\calM$ and the convergence is smooth at any regular point of $\calM$, by the local regularity for Brakke flows, we see that for every $r>1$ there exists $i_0$ such that for all $i\ge i_0$ we have that 
$$ \calM^i_{O,\Lambda(r)} \text{ is } r^{-1} \text{ close in } C^{2,\alpha} \text{ to } \mathcal{C} \text{ on } B_{r}(\bO)\times \left(-r^{2}, - r^{-2}\right).$$

Take now any sequence of regular points $X_i \in \calM^i$ with $X_i \rightarrow (\bO,0)$. From the above it follows that $Z_{\calM^i}(X_i) < \infty$ for $i$ large enough, and there exists a sequence $\{r_i\}_{i=i_0}^\infty$ of positive numbers with $r_i/Z_{\calM^i}(X_i) \to \infty$ such that $\calM^i$ remains $\varepsilon$-cylindrical around $X_i$ at all scales between $Z_{\calM_i}(X_i)$ and $r_i$. 

The rest of the proof in \cite[Proposition 9.1]{CHHW} now works completely analogously: then $\mathcal{M}^i_{X_i, Z(X_i)^{-1}}$ subconverges to an ancient asymptotically cylindrical flow $\mathcal{M}^\infty$ with $Z(O)=1$. By \cite[Theorem 1.5]{CHHW} the limit  $\mathcal{M}^\infty$ is either a round shrinking cylinder, a translating bowl, or an ancient oval. If $\mathcal{M}^\infty$ is the cylinder, then $0$ cannot be
its time of extinction, since $Z(O) = 1$. Therefore, if $\mathcal{M}^\infty$ is either the cylinder or the bowl, it follows that $O$ is a regular point of  $\mathcal{M}^\infty$ and the convergence is smooth. 

The remainder of the argument in the proof of \cite[Proposition 9.1]{CHHW} establishes $(i) - (iv)$. Note that since in the definition of the cylindrical scale we kept track of the axis, we obtain that all the limit flows have the same axis. Note that this remark has been made by the authors already in \cite[p.~21]{CHHW}.

For $(v)$, consider $\beta_0>0$ to be minimal uniform 2-convexity constant of the sphere, the ancient oval, the bowl and the shrinking cylinder. Choose any $0<\beta<\beta_0$ and consider regular points $X_i \in \calM^i$ with $X_i \rightarrow (\bO, 0)$, where $\lambda_1+\lambda_2 \leq \beta H$. Using $(iv)$ we reach a contradiction. 

For $(vi)$, let $\alpha_0$ be the minimum of non-collapsing constants of the shrinking cylinder, the Bowl, the ancient oval and the shrinking round sphere. By $(i)$ all approximating flows have a mean convex neighborhood of $(\bO, 0)$. We can thus for each point $Y=(\by,t)$ on $\calM^i$ consider the quantity $\barr{\alpha}(Y) = r(\by,t) H(\by,t)$, where $r(\by,t)$ is the maximal radius of a ball touching $\calM^i(t)$ at $\by$ from the inside and outside. Choose $0<\alpha<\alpha_0$ and again assume that there are regular points  $X_i \in \calM^i$ with $X_i \rightarrow (\bO, 0)$, where $\barr{\alpha}(X_i) < \alpha$. We reach a contradiction as before.
\end{proof}

\begin{lemma}
\label{lemma-curvature-lower-bound}
There exists a universal constant $\varepsilon_0>0$ such that the following holds. Let $\calM$ be a unit-regular, $n$-integral Brakke flow which has at $O=(\bO,0)$ a nondegenerate neckpinch. If $n=2$ we also assume that $\calM$ is cyclic mod 2. Then for every $\eta>0$ there exists $r>0$ such that 
\begin{equation}\label{eq:lower bound}
R(X)^{-2} \geq \frac{\varepsilon_0}{|t| + \eta |\bx|^2}
\end{equation}
for all regular points $X=(\bx,t) \in \mathcal{M}$ on $B_r(\bO)\times (-r^2,0]$.
\end{lemma}

\begin{proof}
The constant $\varepsilon_0>0$ will be chosen in the course of the proof. We recall the argument in the proof of Proposition \ref{prop-help}. For a given $\eta>0$, we can assume that there is a sequence 
 of regular points $X_i = (\bx_i,t_i) \in \calM$ with $t_i \leq 0$ and $X_i \rightarrow (\bO,0)$, such that
 \begin{equation}\label{eq:lower bound-neg}
R(X_i)^{-2} \leq \frac{\varepsilon_0}{|t_i| + \eta |\bx_i|^2}.
\end{equation}
As in the proof of Proposition \ref{prop-help} it follows that $Z(X_i) < \infty$ for $i$ large enough, and there exists a sequence $\{r_i\}_{i=i_0}^\infty$ of positive numbers with $r_i/Z(X_i) \to \infty$ such that $\calM^i$ remains $\varepsilon$-cylindrical around $X_i$ at all scales between $Z(X_i)$ and $r_i$. 

Then $\mathcal{M}_{X_i, Z(X_i)^{-1}}$ subconverges to an ancient asymptotically cylindrical flow $\mathcal{M}^\infty$ with $Z(O)=1$. Since $\calM$ has a nondegenerate neckpinch at $O$, the limit  $\mathcal{M}^\infty$ is a round shrinking cylinder. Again note that $0$ cannot be
its time of extinction, since $Z(O) = 1$. Even more, the condition $Z(O) = 1$ determines an extinction time $\bar{\tau} > 0$. Note that since the convergence is smooth, there exists  a universal constant $C>0$ such that for $i$ sufficiently large
$$ C^{-1} R(X_i) \leq Z(X_i) \leq C R(X_i)\, .$$
We now choose $\varepsilon_0>0$ such that $C^2 \varepsilon_0 \leq \bar{\tau}/2$. Then \eqref{eq:lower bound-neg} implies that
$$ Z(X_i)^{-2}|t_i| + \eta |Z(X_i)^{-1} \bx_i|^2 \leq C^2\varepsilon_0 \leq \frac{\bar{\tau}}{2} .$$
But this yields that the original singularity of $\calM$ at $O$ converges under this sequence of rescalings to a regular point of $\calM^\infty$, which yields a contradiction for $i$ sufficiently large.
\end{proof}

We now characterise limit flows along sequences at first-time neckpinches.

\begin{lemma}
\label{lemma-limits-uniform}
Let $\calM$ be a unit-regular, $n$-integral Brakke flow which has at $O=(\bO,0)$ a multiplicity one cylinder $\mathcal{C}$ as its tangent flow.   Let $\calM^i$ be a sequence of unit-regular, $n$-integral Brakke flows converging to $\calM$ in a neighborhood of $O$,
having as tangent flows at singular points  $\bar X_i=({\bf \bar x_i}, \bar t_i) \rightarrow O$ multiplicity one cylinders $\mathcal{C}_i$ {\rm (}with possibly different axes{\rm)}. If $n=2$ we also assume that all flows are cyclic mod 2.
Assume that there is a neighborhood $U$ of $O$ such that $\bar{t}_i$ is the first singular time of $\calM^i$ in $U$ for all $i$ sufficiently large. Then for every $\eta >0$ there exist an $r_0 > 0$ and $i_0$ sufficiently large, such that for $i\geq i_0$, every $X = (\bx,t)$ on $\calM_i$ with $t < \bar{t}_i$ and $X \in B_{r_0}(\bO) \times (-r_0^2, r_0^2)$, the rescaled flow $ \calM^i_{X, R(X)^{-1}}$ is smoothly $\eta$-close to either a shrinking cylinder or a translating bowl, where the axis is aligned with the axis of the neckpinch of $\calM$ at $O$.
\end{lemma}

\begin{proof}
From the proof of Proposition \ref{prop-help} we have that for every $r>1$ there exist $\Lambda(r)>0$ and $i_0$ such that for all $i\ge i_0$ we have that 
\begin{equation}\label{eq:closeness_cylinder}
 \calM^i_{O,\Lambda(r)} \text{ is } r^{-1} \text{ close in } C^{2,\alpha} \text{ to } \mathcal{C} \text{ on } B_{r}(\bO)\times \left(-r^{-2}, - r^{2}\right).
\end{equation}

Let us denote with $a \in \mathbb{RP}^{n+1}$ the axis of $\mathcal{C}$ and with $a_i \in \mathbb{RP}^{n+1}$ the axes of the cylindrical tangent flows $\mathcal{C}_i$ of $\calM^i$ at $\bar X_i$. For every $\varepsilon >0$ there exist $i_0$ such that $\bar X_i \in B_\varepsilon(\bO)\times (-\varepsilon^2,\varepsilon^2)$ for $i\geq i_0$. We now want to argue that $a_i \rightarrow a$.

We can assume that for a given $r>0$ and $i \geq i_0$ the flows $\calM^i_{O,\Lambda(r)}$ satisfy \eqref{eq:closeness_cylinder}. Note further that we can assume that the area ratios of $\calM_i$ are uniformly bounded. Thus by the {\L}ojasiewicz-Simon inequality of Colding-Minicozzi \cite{CM} we see that the flows $\calM^i_{O,\Lambda(r)}$ are $(r/2)^{-1}$ close in  $C^{2,\alpha}$-norm to $\mathcal{C}_i$ on $B_{r/2}(\bO)\times \left(-(r/2)^2, -(r/2)^{-2} \right)$. This implies that $a_i \rightarrow a$.

Note furthermore that the {\L}ojasiewicz-Simon inequality of Col\-ding-\-Mini\-cozzi \cite{CM} implies a uniform decay rate of $\calM^i_{O,\Lambda(r)}$ to each of its cylindrical tangent flows at $\bar{X}_i = (\bar{\bx}_i,\bar{t}_i)$ in the time interval $\big(\bar{t}_i, - (r/3)^2\big)$.
We can thus consider the flows $\tilde{\calM}^i$ obtained from $\calM^i_{O,\Lambda(r)}$ by translating by $-\bar X_i$ and rotating such that $a_i =a$ for all $i>i_0$. Note that this family of flows decays uniformly towards their tangent flow $\mathcal{C}$ for $t\in (-1,0)$. 

Then as in the proof of $(iv)$ in Proposition \ref{prop-help}  (compare it with the proof of  \cite[Proposition 9.1]{CHHW}), we obtain that for any sequence of points $X_i= (\bx_i,t_i)$ on $\tilde{\calM}^{j(i)}$ with $t_i < 0$ such that 
$X_i \rightarrow O$,  a sequence of rescalings $\tilde{\calM}^{j(i)}_{X_i, R^{-1}(X_i)}$ converges to either a round shrinking sphere, a round shrinking cylinder, a translating bowl, or an ancient oval, with the axis aligned with $a$. But note that the limit can neither be the shrinking sphere or the ancient oval, since all the flows $\tilde{\calM}_{j(i)}$ have a first time cylindrical singularity at $O$. Thus the limit can only be the translating bowl or the shrinking cylinder with the axis aligned with $a$.
\end{proof}

\begin{lemma}
\label{lemma-all-cylinder}
Assume $\{M_t\}_{t\in [0,T)}$ is a smooth mean curvature flow with a nondegenerate neckpinch at $(\bO,T)$, which is an isolated singularity, such that the axis of the neckpinch is the $\bx_1$-axis. Then there exists $r>0$ such that on 
$$C_r=\left \{\bx \in \R^{n+1}\, \Big|\,  |\bx_1|\leq r \text, \sum\nolimits_{i=2}^{n+1} \bx_i^2 \leq r^2\right\}$$ 
and all $t \in (T-r^2,T)$, $M_t$ can be written as a cylindrical graph over the $\bx_1$-axis and the heightfunction $u = u(\bx_1, \omega)$ satisfies $u = O(\sqrt{T-t})+ o(|\bx_1|)$. Furthermore, on $C_r\setminus \{\bO\}$ the flow converges smoothly to a limiting hypersurface $M_T$.  $M_T$ can be written as a cylindrical graph over the $\bx_1$-axis for all $x_1 \in (-r,r)\setminus \{0\}$ and the heightfunction $u = u(\bx_1, \omega)$, with $\omega \in \mathbb{S}_1^n$, satisfies $u = o(|\bx_1|)$.
\end{lemma}
\begin{proof} 
Since $(\bO,T)$ is a nondegenerate neckpinch, together with Lemma \ref{lemma-curvature-lower-bound}, we see that for every $\eta > 0$ there exists an $r >0$ sufficiently small, so that around every $X= (\bx,t) \in \calM \cap (C_r \times (T-r^2,T))$ we can assume that $M_t$ is close to a shrinking cylinder around the $\bx_1$-axis with radius less than $C\sqrt{T-t+\eta |\bx_1|^2}$, for some uniform constant $C$. Indeed, to see that, let us argue by contradiction. Assume there exist $\eta > 0$ and a sequence $X_i = (\bx_i,t_i) \to ({\bf 0},T)$, so that for any $i$, $M_{t_i}$ around $X_i$ is not close to a cylinder around the $\bx_1$ axis with radius less than $C_0\, \sqrt{T-t_i+\eta|{\bf x}_1^i|^2}$, where $\bx_i = (\bx^i_1, \dots \bx^i_n)$ and $C_0 = \epsilon_0^{-1/2}$ with $\epsilon_0$ is as in Lemma \ref{lemma-curvature-lower-bound}. Rescale $\calM$ around $X_i$ by $R(X_i)^{-1}$. By Lemma \ref{lemma-curvature-lower-bound} we have that the $\lim_{i\to\infty} R(X_i)^{-1} = \infty$. Since $\calM$ has a nondegenerate neckpinch at $({\bf 0},T)$, by Definition \ref{def-nondegenerate} we have that the sequence of rescaled flows $\calM_{X_i,R(X_i)^{-1}}$ converges to a round cylinder. This in particular means that   $M_{t_i}$, around $X_i$, is close to a cylinder around the $\bx_1$-axis with radius less than $R(X_i) \le C_0\sqrt{T-t_i+\eta |\bx^i_1|^2}$. Hence, we get contradiction.

Since the tangent flow of $\calM$ at $(\bO,T)$ is a shrinking cylinder, we can furthermore assume that on $\{(\bx, t)\, | \, T-t \geq \eta \,\bx^2 \} \cap C_r$ the flow $\calM$ is close to a shrinking cylinder with unit multiplicity. A direct continuation argument gives that $M_t\cap C_r$ extends out to $\{|\bx_1| =r\}$ as a smooth cylindrical graph. Since this is true for all $t<T$, together with the assumption that $(\bO,T)$ is an isolated singularity, this yields that the same extends to $t=T$.
\end{proof}

A direct consequence of Lemma  \ref{lemma-all-cylinder} is Theorem \ref{thm-no-vanishing}.

\begin{proof}[Proof of Theorem \ref{thm-no-vanishing}]
Without loss of any generality assume there is only one isolated nondegenerate neckpinch singularity at $({\bf 0},T)$. Then by Lemma \ref{lemma-all-cylinder} there exists an $r > 0$ so that on $C_r\backslash \{{\bf 0\}}$ the flow converges smoothly to $M_T$, and $M_T$ can be written as a cylindrical graph over the $\bx_1$ axis for all $|\bx_1| <r$. This immediately concludes the Theorem.
\end{proof}

\section{Stability of neckpinch singularities}
\label{sec-neck}

In this section we first prove a local version of Theorem \ref{thm-main1}, which roughly says that if under certain assumptions a tangent flow at the first singular time of a smooth flow $\{M_t\}$ is a round cylinder, the same is true for a mean curvature flow starting sufficiently close by. Theorem \ref{thm-main1} then follows as a direct corollary. This tells us neckpinch type singularities are dynamically stable.

\begin{theorem}
\label{thm-main1-local}
Let $\calM$ be a unit regular Brakke flow that has an isolated singularity at $O=(\bO,0)$, which is a neckpinch.  Assume that the neckpinch locally disconnects the manifold into two pieces none of which disappears at time $0$.
Let $\calM^i$ be a sequence of unit regular Brakke flows converging to $\calM$. If $n=2$ we also assume that all flows are cyclic mod 2. Then there exits $r>0$ and $i_0 \in \mathbb{N}$ such that for $i \geq i_0$, the first singularity of $\calM^i$ in $B_{r}(\bO)\times (-r^2,r^2)$ occurs in $B_{r/2}(\bO)\times (-r^2/4,r^2/4)$ and is a neckpinch.
\end{theorem}

\begin{remark}
 Note that the proof also yields that the first singularity of $\calM^i$ in $B_{r}(\bO)\times (-r^2,r^2)$ occurs in $B_{r/2}(\bO)\times (-r^2/4,r^2/4)$, is a neckpinch and locally disconnects the manifold into two pieces none of which disappears at the singular time. But it appears to be a subtle question to show that this singularity is also isolated.
\end{remark}

\begin{proof}
We can assume that the axis of the neckpinch of $\calM$ at $O$ is the $\bx_1$-axis. Let  $r_0=r$ and $i_0$ be as in Proposition \ref{prop-help}. Since by our assumption the neckpinch locally disconnects $\calM(0)$ into two pieces none of which disappears at time $0$, we can choose $0<4r_1<r_0$ such that $\calM$ is smooth on $(B_{4r_1}(\bO) \setminus B_{r_1}(\bO)) \times [-16r_1^2,0]$, and consists of two smooth  connected components, one in the halfspace $\{\bx_1 <0\}$ and one in the halfspace $\{\bx_1 > 0\}$, and both having boundary on $\partial B_{4r_1}(\bO)$ and $\partial B_{r_1}(\bO)$. Since the tangent flow of $\mathcal{M}$ at $({\bf 0},0)$ is a round multiplicity one cylinder, by taking $r_1$ sufficiently small we can furthermore assume that $\calM$ is smoothly close to the shrinking cylinder $\calC$ on $B_{4r_1}(\bO)\times [-16 r_1^2, - 4 r_1^2]$.

Note that since $\calM$ is unit regular, there exists $\eta>0$ such that $\calM$ is smooth and non-vanishing on $(B_{4r_1}(\bO) \setminus B_{r_1}(\bO)) \times [-16r_1^2, \eta)$. 

Let $\Gamma$ be a scaling of the Angenent torus, rotationally symmetric with respect to the $\bx_1$-axis, and symmetric with respect to the $\{\bx_1=0\}$-plane such that it encloses the cylinder $\mathbb{S}^{n-1}(\sqrt{8(n-1)})$. Note that this is the cylinder with double the radius than the standard shrinking cylinder at time $-1$. Let $\mathcal{G} = (\Gamma_t)_{t \in [-1,\gamma)}$ be the maximal evolution of $\Gamma$ starting at time $-1$. Note that since $\Gamma$ is disjoint from the shrinking cylinder at time $-1$, we have that $\gamma >0$.

Choose $\lambda_0>0$ such that $\lambda_0^2\eta/2 \geq 2 \gamma$ and $\lambda_0r_1 \geq 2(1+\sqrt{\gamma})$. Let $R_0$ be such that $\Gamma \subset B_{R_0}(\bO)$. Since the singularity of $\calM$ at $O$ is a neckpinch, for every $\epsilon > 0$ there exists $\lambda \geq \lambda_0$ such that $\calM_{O,\lambda}$ is $\varepsilon$-close in $C^2$-norm to the shrinking cylinder on $B_{100R_0}(\bO)\times [-4, -1/4]$. We can choose $\varepsilon>0$ such that any surface which is $2\varepsilon$-close in $C^2$-norm to the shrinking cylinder at time $t=-1$ is still enclosed by $\Gamma$. 

We can now choose $i_1\geq i_0$ such that for all $i\geq i_1$ we have that $\calM^i_{O,\lambda}$ is $2\varepsilon$-close in $C^2$ to the shrinking cylinder on $B_{100R_0}(\bO)\times [-4, -1/4]$ and that $\calM^i$ is $2\varepsilon$-close in $C^2$ to $\calM$ on $(B_{3r_1}(\bO) \setminus B_{2r_1}(\bO)) \times [-8r_1^2,\eta/2]$ as well as $2\varepsilon$-close in $C^2$  to  $\calM$ on $B_{3r_1}(\bO)\times [-8 r_1^2, - 2/\lambda^2]$. 

Note that these choices imply that $\calM^i_{O,\lambda}$ is enclosed by $\Gamma$ at time $t=-1$. Scaling back we see that $\calM^i(\lambda^{-2})$ is enclosed by $\lambda \Gamma$, and the evolution of $\lambda \Gamma$ contracts to $\bO$ at time $\lambda^{-2}\gamma$. We now aim to show that $\calM^i$ has a neckpinch singularity before $\lambda^{-2}\gamma$ in $B_{2r_1}(\bO) \times [- 4 r_1^2, \lambda^{-2}\gamma]$.

We can first note that the tangent flow at the first singularity (if it occurs) of $\calM^i$ in $B_{2r_1}(\bO) \times [- 4 r_1^2, \lambda^{-2}\gamma]$ has to be either a sphere or a cylinder due to the uniform 2-convexity (by part (v) of Proposition \ref{prop-help}). Since $\calM^i$ is smoothly close to $\calM$ on $(B_{3r_1}(\bO) \setminus B_{2r_1}(\bO)) \times [-8r_1^2,\eta/2]$ and on 
$B_{3r_1}(\bO)\times [-8 r_1^2, - 2/\lambda^2]$, it can't be a sphere. 

We now show that a singularity has to occur in $B_{2r_1}(\bO) \times [- 4 r_1^2, \lambda^{-2}\gamma]$. Choose points $z_1 \in \partial B_{3r_1}(\bO) \cap \{\bx_1<0\}$ and $z_2 \in \partial B_{3r_1}(\bO) \cap \{\bx_1>0\}$ both enclosed by $\calM^i(\eta/2)$. Note that since $\calM^i$ is mean convex in $B_{r_0}(\bO)\times 
(-r_0^2,r_0^2)$ we have that $z_1,z_2$ are also enclosed by $\calM^i(-4r_1^2)$. Note that since $\calM^i(-4r_1^2)$ is smoothly close to a shrinking cylinder at time $-4r_1^2$, both $z_1,z_2$ are in the same connected component of $B_{3r_1}(\bO) \setminus \calM^i(-4r_1^2)$. 

\begin{claim}
\label{claim-components}
 The points $z_1,z_2$ are in different connected components of $B_{3r_1}(\bO) \setminus \calM^i(\eta/2)$ \end{claim}
 
 \begin{proof}
Assume there would be a curve $\alpha$ connecting $z_1,z_2$ in $B_{3r_1}(\bO) \setminus \calM^i(\eta/2)$. Note that  since $\calM^i$ is mean convex in $B_{r_0}(\bO)\times (-r_0^2,r_0^2)$, $\alpha$ has to also connect $z_1,z_2$ in $B_{3r_1}(\bO) \setminus \calM^i(t)$ for $t\in [-4r_1^2,\eta/2]$ .  Recall that $\calM^i(\lambda^{-2})$ is enclosed by $\lambda \Gamma$, and $\alpha$ has to run 'through' the Angenent torus $\lambda \Gamma$ at time $-\lambda^{-2}$. This yields a contradiction since the $\lambda \Gamma$ shrinks away at time $\lambda^{-2} \gamma$ and $\eta/2 > \lambda^{-2} \gamma$.
\end{proof}
To complete the statement of the theorem we can choose $r = 4r_1$ and relabel $i_1$ with $i_0$.
\end{proof}

\begin{proof}[Proof of Theorem \ref{thm-main1}]
Theorem \ref{thm-main1} is a direct consequence of Theorem \ref{thm-main1-local}. Assume without loss of any generality that $(M_t)_{t \in [0,T)}$ has an isolated singularity at $O =(\bO,0)$. Let $(M^i_t)_{t \in [0,T_i)}$ be a sequence of smooth flows such that $M^i_0 \rightarrow M_0$ in $C^2$. By elliptic regularization, \cite{I94}, there exist integral Brakke flows $\calM, \calM^i$ which coincide with the smooth flows $(M_t)_{t \in [0,T)}$ and $(M^i_t)_{t \in [0,T_i)}$, respectively, as long as they exist. Note that the flows constructed via elliptic regularization are always unit regular, see \cite{SW}, and for $n=2$ we can also assume that they are cyclic mod 2. Furthermore, since $(M_t)_{t \in [0,T)}$ only has neckpinches at time $T$, by \cite{HW}, there exists $\delta>0$ such that weak evolution $\calM$ is unique for $t\in [0,T+\delta]$. But this implies that $\calM^i$ converge weakly to $\calM$ in a neighborhood of $O$. Then the statement follows from Theorem \ref{thm-main1-local}.
\end{proof}

\section{Stability of nondegenerate neckpinches }
\label{sec-nondeg-neck}

In this section we prove Theorem \ref{thm-main2}. It roughly says that if the mean curvature flow starting at $M_0$ develops a nondegenerate neckpinch at the first singular time, then the mean curvature flow starting at any sufficiently small perturbation of $M_0$ will also develop only nondegenerate neckpinches at the first singular time. 

We will first state a local version of Theorem \ref{thm-main2}.

\begin{theorem}\label{thm-nondegenerate-stability-local}
Let $\calM$ be a unit regular $n$-integral Brakke flow which has an isolated singularity at $O=(\bO,0)$ which is a nondegenerate neckpinch. Let $\calM^i$ be a sequence of unit regular Brakke flows converging to $\calM$ in a neighborhood of $O$. If $n=2$ we also assume that all flows are cyclic mod 2. Then there exits $r>0$ and $i_0 \in \mathbb{N}$ such that for $i \geq i_0$, the first singularity of $\calM^i$ in $B_{r}(\bO)\times (-r^2,r^2)$ occurs in $B_{r/2}(\bO)\times (-(r/2)^2,(r/2)^2)$ and is a nondegenerate neckpinch.
\end{theorem}

\begin{proof}
 We first note that by Lemma \ref{lemma-all-cylinder} there is an $r_1>0$ such that in $B_{r_1}(\bO)\times (-r^2_1,0]$, the flow $\calM$ can be written as smooth cylindrical graph over the axis $a$ of the neckpinch at $O$ and the flow disconnects into two pieces at time $0$, non of which disappears. Thus by Theorem \ref{thm-main1-local} there exists $0<r_2<r_1$ and $i_2>0$ such that for $i>i_2$ the first singularity of $\calM^i$ in $B_{r_2}(\bO)\times (-r_2^2,r_2^2)$ occurs  at $X_i = (x_i,t_i)  \in B_{r_2/2}(\bO)\times (-(r_2/2)^2,(r_2/2)^2)$   and is a neckpinch. 
For $i>i_2$ assume that a special limit flow of $\calM_i$ at $X_i$ is a translating bowl. We can assume that $r_2$ is smaller than $r_0$ given in  Lemma \ref{lemma-limits-uniform}. Applying Lemma \ref{lemma-limits-uniform} together with the smoothness of the flow $\calM^i$ before $t_i$ we see that there has to be a continuous time-directed curve of points $\gamma(t) $ on $\calM^i(t)$ 
 where the normal vector 
of $\calM^i(t)$ at $\gamma(t)$ is parallel to the axis $a$ (these points correspond to the tip of the limiting  bowl soliton). 

Even more, by Lemma  \ref{lemma-limits-uniform} and a continuity argument, the curve $\gamma(t)$ can be extended backward in time until it hits the parabolic boundary of $B_{r_2}(\bO)\times (-r_2^2,0)$. On the other hand, we  can assume that in  a neighborhood of the parabolic boundary of $B_{r_2}(\bO)\times (-r_2^2,0)$ the flow $\calM^i$ is smoothly close to $\calM$ which  is cylindrical there by Lemma \ref{lemma-all-cylinder}. This yields   a contradiction.
\end{proof}

\begin{proof}[Proof of Theorem \ref{thm-main2}]
This is a direct application of Theorem \ref{thm-nondegenerate-stability-local} as in the proof of Theorem \ref{thm-main1}
\end{proof}

\section{Local Type I and Type II  singularities}\label{sec-typeI}

In this section we give a proof of Theorem \ref{prop-typeI}. This gives us criteria which guarantee that a neckpinch singularity is either Type I or Type II.

\begin{proof}[Proof of Theorem \ref{prop-typeI}] The claim that a nondegenerate neckpinch is always Type I easily follows from the following fact. Assume that $\calM$ has a nondegenerate neckpinch singularity at $X = (\bx_0,t_0)$. We argue by contradiction. If the singularity at $X$ were Type II, then by Hamilton's procedure we would be able to carefully choose a sequence of points $\bx_i\in M$, where $\bx_i\to \bx_0$, and a sequence of times $t_i\to t_0, t_i<t_0$, so that a sequence of rescaled solutions  around $(\bx_i,t_i)$ by  $|A|(\bx_i,t_i)$, smoothly converges to a blow up limit which is a non-trivial eternal solution. This would contradict our assumption that we have a nondegenerate neckpinch at $X$ and hence every limit flow must be a round cylinder. 

Assume that $\calM$ has a degenerate neckpinch singularity at $X = (\bx_0,t_0)$, but satisfies a local Type I bound, i.e. there exists $\delta>0$ and $C>0$ such that
$$|A|^2(\mathbf{x},t) \leq \frac{C}{t_0-t}$$
for all $(\bx,t) \in B_\delta(\mathbf{x}_0)\times (t_0-\delta^2,t_0)$. Pick a sequence of points $X_i=(\bx_i,t_i) \rightarrow X$ with $t_i\leq t_0$, $t_i < t_0$ if $\bx_i \neq \bx_0$ and scaling factors $\lambda_i \rightarrow \infty$. 
Note that due to the Type I assumption, we have that the second fundamental form of rescaled flow $\calM_{X_i,\lambda_i}$, call it  $A_i$, satisfies the bound
\begin{equation}
\label{eq-norm}
|A_i| \le \frac{C}{\lambda_i(t_0-t_i)-t},
\end{equation}
implying that the flow $\calM_{X_i,\lambda_i}$ always smoothly converges to an ancient flow $\calM'$, as $i\to\infty$, for all times $t < \liminf_{i\to \infty} \lambda_i(t_0-t_i)$. By parts (v) and (vi) in Proposition \ref{prop-help} we know that $\calM'$ is uniformly 2-convex and $\alpha$-noncollapsed, for some $\alpha > 0$. By \cite{BC} we know  that if $\calM'$ is nontrivial, it is either a plane, a shrinking sphere, a shrinking cylinder or the translating bowl. If
$\calM'$ were  the shrinking sphere, $\calM$ would have to be strictly convex for all times sufficiently close to a singular time $t_0$. That would imply the tangent flow at $X$ would have to be of strictly convex type, and thus could not be a neckpinch, hence contradiction.

So the only thing we need to rule out is the translating bowl. Up to subsequences we can distinguish the following three cases:
\begin{itemize}
\item[(1)] $\lim_{i \rightarrow \infty} \lambda_i^{-1} |A|(\bx_i,t_i) = 0$\,,
\item[(2)] $\lim_{i \rightarrow \infty} \lambda_i^{-1} |A|(\bx_i,t_i) = \Lambda >0$\,,
\item[(3)] $\lim_{i \rightarrow \infty} \lambda_i^{-1} |A|(\bx_i,t_i) = + \infty$\,.
\end{itemize} 
Note that in case (1), the limit flow is trivial. 

In cases (2) or (3), by \eqref{eq-norm} we have
\begin{equation}
\label{eq-est-A}
\lambda_i^{-1}|A|(\bx_i,t_i) \le \frac{C}{\lambda_i(t_0-t_i)},
\end{equation}
implying that the $\lim_{i\to \infty} \lambda_i(t_0-t_i) \le C/\Lambda =: T \ge 0$. Hence, the limit flow is an ancient flow and it satisfies a type I estimate of the form
 $$|A|^2(\mathbf{x},t) \leq \frac{C}{T-t}$$
 for all $t<T$. This immediately rules out the translating bowl.
 \end{proof}

\bibliographystyle{alpha}
\bibliography{stability-neckpinch}

\end{document}